\documentclass[graybox]{svmult}

\usepackage{mathptmx}       
\usepackage{helvet}         
\usepackage{courier}        
\usepackage{type1cm}        
\usepackage[bottom]{footmisc}

\usepackage{amsfonts}
\usepackage{amsmath}
\usepackage{amssymb}
\usepackage{array}
\usepackage{stmaryrd}
\usepackage{graphicx}

\usepackage{mathabx} 
\usepackage{bbm}

\usepackage{tikz}
    \usepackage{color}  
     \definecolor{red}{rgb}{0.9,0,0}
     \definecolor{green}{rgb}{0,0.6,0}
     \definecolor{rb}{rgb}{0.6,0,0.2}     
     \definecolor{blue}{rgb}{0,0,0.7}
\renewcommand{\sim}{\simeq}
\renewcommand{\epsilon}{\varepsilon}

\newcommand{\pt}{\partial}
\newcommand {\eps} {\varepsilon}

\newcommand{\bmu}{\boldsymbol{\mu}}
\newcommand{\bnu}{{\boldsymbol{\nu}}}

\newcommand{\vvvertN}{{/\hspace{-2.9pt}/\hspace{-2.9pt}/\!}}

\newcommand{\proofend}{\hfill$\Box$}

\newcommand {\beq} {\begin{equation}}
\newcommand {\eeq} {\end{equation}}
\newcommand {\beqa} {\begin{eqnarray}}
\newcommand {\eeqa} {\end{eqnarray}}
\newcommand {\beqann} {\begin{eqnarray*}}
\newcommand {\eeqann} {\end{eqnarray*}}

\numberwithin{equation}{section}
\renewcommand{\theequation}{\arabic{section}.\arabic{equation}}
\numberwithin{remark}{section}
\numberwithin{lemma}{section}
\newtheorem{theorem_}[lemma]{Theorem}
\newtheorem{corollary_}[lemma]{Corollary}

\begin{document}

\title*{Improved energy-norm a posteriori error estimates for singularly perturbed reaction-diffusion problems\\ on
anisotropic meshes}
\titlerunning{Improved energy-norm a  posteriori error estimates on
anisotropic meshes}
\author{Natalia Kopteva}
\institute{Natalia Kopteva \at University of Limerick, Limerick, Ireland; \email{natalia.kopteva@ul.ie}
}
%
%
\maketitle

\abstract*{In the recent article \cite{Kopt_17_NM}
the author obtained
residual-type a posteriori error estimates in the energy norm
for singularly perturbed semilinear reaction-diffusion equations
on anisotropic triangulations.
The error constants in these estimates are
independent  of the diameters and the aspect ratios of mesh elements
and of the small perturbation parameter.
The purpose of this note is to improve the weights in the jump residual part of the estimator.
This is attained by using a novel sharper version of the scaled trace theorem for anisotropic elements, in which the hat basis functions are involved as weights.}

\abstract{In the recent article \cite{Kopt_17_NM}
the author obtained
residual-type a posteriori error estimates in the energy norm
for singularly perturbed semilinear reaction-diffusion equations
on unstructured anisotropic triangulations.
The error constants in these estimates are
independent  of the diameters and the aspect ratios of mesh elements
and of the small perturbation parameter.
%
The purpose of this note is to improve the weights in the jump residual part of the estimator.
This is attained by using a novel sharper version of the scaled trace theorem for anisotropic elements, in which the hat basis functions are involved as weights.}

\renewcommand{\theequation}{\arabic{section}.\arabic{equation}}

\section{Introduction}
\label{sec:1}
Consider finite element approximations to
singularly perturbed semilinear reaction-diffusion equations of the form
\begin{align}
Lu:= -\epsilon^2 \triangle u + f(x,y;u) = 0\quad\mbox{for}\;\;(x,y)\in \Omega, \qquad u=0\quad \mbox{on}\;\; \partial \Omega,
\label{eq1-1}
\end{align}
posed in a,
possibly non-Lipschitz,
polygonal domain $\Omega\subset\mathbb{R}^2$.
Here $0<\epsilon \le 1$. We also assume that
$f$ is continuous {on $\Omega \times \mathbb{R}$ and satisfies $f(\cdot; s) \in L_\infty(\Omega)$ for all $s \in \mathbb{R}$},
and
the one-sided Lipschitz condition $f(x,y; v)-f(x,y;w) \ge C_f [v-w]$ whenever $v\ge w$,
with some constant $C_f \ge 0$.
Then there is a unique solution
$u\in W_\ell^2(\Omega)\subseteq W_q^1\subset C(\bar\Omega)$ for some  $\ell>1$ and $q>2$
\cite[Lemma~1]{DK14}.
We additionally assume that $C_f+\varepsilon^2\ge 1$
(as \eqref{eq1-1} can always be reduced to this case by a division by $C_f+\varepsilon^2$).

For this problem, the recent articles \cite{Kopt15,Kopt_17_NM} gave
residual-type
a posteriori error estimates
on unstructured anisotropic meshes.
In particular, in \cite{Kopt_17_NM} the error was estimated in the energy norm $\vvvert \cdot \vvvert_{\eps\,;\Omega}$,
which is an appropriately scaled $W^{1}_2(\Omega)$ norm naturally associated with our problem, defined for any ${\mathcal D}\subseteq \Omega$ by
$
\vvvert v \vvvert_{\eps\,;{\mathcal D}}:=
\Bigl\{\eps^2\|\nabla v \|^2_{2\,;{\mathcal D}}
+\|v \|^2_{2\,;{\mathcal D}} \Bigr\}^{1/2}
$.
%
Linear finite elements were used to discretize \eqref{eq1-1}
with
a piecewise-linear finite element space
 $S_h \subset H_0^1(\Omega)\cap C(\bar\Omega)$
relative to a triangulation $\mathcal T$, and the
the computed solution $u_h \in S_h$ satisfying
\begin{align}
\eps^2 \langle\nabla u_h,\nabla v_h\rangle+\langle f_h^I,v_h\rangle=0\quad
\forall\;v_h \in S_h,
\qquad\quad
f_h(\cdot):=f(\cdot;u_h).
\label{eq1-2}
\end{align}
Here $\langle\cdot,\cdot\rangle$ denotes the $L_2(\Omega)$ inner product,
and $f_h^I$ is the standard piecewise-linear Lagrange interpolant of $f_h$.

To give a flavour of the results in \cite{Kopt_17_NM},
assuming that all mesh elements are anisotropic, 
 one estimator reduces to
\begin{align}
\notag
\vvvert u_h-u \vvvert_{\eps\,;\Omega}\le C
\,\,
\Bigl\{&\sum_{z\in{\mathcal N}}\min\{|\omega_z|,\,\lambda_z\}
\,\bigl\|\eps \llbracket\nabla u_h\rrbracket \bigr\|^2_{\infty\,;\gamma_z}
\\\label{first result}
{}+{}&
\sum_{z\in{\mathcal N}}
\bigl\|\min\{1,\,H_z\eps^{-1}\}\,f_h^I\bigr\|^2_{2\,;\omega_z}
+\,\bigl\|f_h-f_h^I\bigr\|^2_{2\,;\Omega}\Bigr\}^{1/2}\!,
\end{align}
where $C$ is 
{independent   of the diameters and the aspect ratios} of elements in $\mathcal T$,
and of $\eps$.
Here
$\mathcal N$ is the set of nodes 
in $\mathcal T$,
$\llbracket\nabla u_h\rrbracket$ is the standard jump in the normal derivative of $u_h$ across an element edge,
$\omega_z$ is the patch of elements surrounding any $z\in\mathcal N$,
$\gamma_z$ is the set of edges in the interior of $\omega_z$, $H_z={\rm diam}(\omega_z)$,
and $h_z\sim H_z^{-1}|\omega_z|$.

A version of \eqref{first result} obtained in \cite{Kopt_17_NM} involves a somewhat surprising weight
$\lambda_z=\eps H_z^2 h_z^{-1}$
at the jump residual terms. The main purpose of this note is to improve the jump residual part of the latter estimator and establish its sharper version with a more natural
$\lambda_z=\eps H_z$.
This will be attained by employing a novel sharper version of the scaled trace theorem,
 in which the hat basis functions are involved as weights (see 
 Remark~\ref{rem_sec_trace}).
As the improvement that we present here applies to the jump residual terms only, we restrict our analysis to these terms.

Note that the new shaper version of the jump residual part of the estimator works not only for \eqref{first result} (see Theorem~\ref{lem_J_gen} below),
but can be also combined with
a shaper bound for the interior residual terms given by \cite[Theorem~6.2]{Kopt_17_NM}. The latter
is more intricate and was obtained under some additional assumptions on the mesh, so we shall not
give it here.
Comparing it to \eqref{first result}, roughly speaking, the weight $\min\{1,\,H_z\eps^{-1}\}$ is replaced by a sharper
$\min\{1,\,h_z\eps^{-1}\}$
with a few  additional terms included.

Note also that a similar improved jump residual part of the estimator is also obtained in
\cite[(1.2)]{Kopt17} using an entirely different (and more complicated in the context of residual-type estimation) approach for a version of \eqref{eq1-2}
(with a special anisotropic quadrature used for the reaction term).

Our interest in locally anisotropic meshes is due to that they
offer an efficient way of computing
reliable numerical approximations
of layer solutions.
(In the context of
\eqref{eq1-1} with $\eps\ll 1$,
see, e.g.,
\cite{Kopt_MC_07,Kopt_EOR,RStTob} and references therein.)
But such anisotropic meshes are frequently constructed a priori or by heuristic methods, while the majority of
available a posteriori error estimators assume shape regularity of the mesh
\cite{AinsOd_2000}.
In the case of shape-regular triangulations, residual-type a posteriori error estimates for equations of type~\eqref{eq1-1}
were proved
in \cite{Ver98c} in the energy norm, and more recently
in \cite{DK14} in the maximum norm.
The case of anisotropic meshes having a tensor-product structure was addressed in \cite{Siebert96}
 for the Laplace equation 
and in \cite{Kopt08,ChKopt}
for problems of type \eqref{eq1-1}, 
with the error  estimators given, respectively, in the $H^1$ norm and the maximum norm.
For unstructured anisotropic meshes, a posteriori error estimates can be found
 in \cite{Kunert2000,KunVer00}
for the Laplace equation in the $H^1$ norm, and in  \cite{Kunert2001,KunVer00}
for a linear constant-coefficient version of \eqref{eq1-1}
in the energy norm.

Note that the error constants in the  estimators of \cite{Kunert2000,Kunert2001,KunVer00} involve the so-called matching functions; the latter depend on the unknown error
and take moderate values only when the grid is either isotropic, or, being anisotropic, is aligned correctly to the solution,
while, in general, they may be as large as mesh aspect ratios.
The presence of such matching functions in the estimator is clearly undesirable. It is entirely avoided in the more recent papers \cite{Kopt15,Kopt_17_NM,Kopt17}, as well as here.

The paper is organized as follows.
In \S\ref{sec_nodes} and \S\ref{sec_trace}, we respectively describe our triangulation assumptions
and
give a novel shaper version of  the scaled trace theorem for anisotropic elements. In \S\ref{sec_error_an},
we dervie the main result of the paper, a new shaper  jump residual part of the estimator.
A simplified version of this analysis is given in \S\ref{sec_struct} for partially structured anisotropic
meshes, while more general anisotropic meshes are addressed in \S\ref{ssec_gen_anal}.
\smallskip

{\it Notation.}
We write
 $a\sim b$ when $a \lesssim b$ and $a \gtrsim b$, and
$a \lesssim b$ when $a \le Cb$ with a generic constant $C$ depending on $\Omega$ and
$f$,
but 
%
not
 on either $\eps$ or
 the diameters and the aspect ratios of elements in $\mathcal T$.
  Also, for $\mathcal{D}\subset\bar\Omega $, $1 \le p \le \infty$, and $k \ge 0$,
  let $\|\cdot\|_{p\,;\mathcal{D}}=\|\cdot\|_{L_p(\mathcal{D})}$
  and $|\cdot|_{k,p\,;\mathcal{D}}=|\cdot|_{W_p^k(\mathcal{D})}$, where $|\cdot |_{W_p^k(\mathcal{D})}$
  is the standard Sobolev seminorm, 
  and ${\rm osc}(v;\mathcal{D})=\sup_{\mathcal{D}}v-\inf_{\mathcal{D}}v$ for $v\in L_\infty(\mathcal{D})$.

\section{Basic triangulation assumptions}\label{sec_nodes}
We shall use $z=(x_z,y_z)$, $S$ and $T$ to respectively denote particular mesh nodes, edges and elements,
while $\mathcal N$, $\mathcal S$ and $\mathcal T$ will respectively denote their sets.
For each $T\in\mathcal T$,
let
$H_T$ be the maximum edge length 
and $h_T:= 2 H_T^{-1}|T|$ be the minimum height in $T$.
For each $z\in\mathcal N$, let
$\omega_z$ be the patch of elements surrounding any $z\in\mathcal N$,
${\mathcal S}_z$ the set of edges originating at $z$,
and
\beq\label{ring_gamma}
H_z:={\rm diam}(\omega_z),\quad h_z:=H_z^{-1}|\omega_z|,
\quad
\gamma_z:={\mathcal S}_z\setminus\pt\Omega,\quad
\mathring{\gamma}_z:=\{S\subset\gamma_z: |S|\lesssim h_z\}.
\eeq
Throughout the paper we make the following triangulation assumptions.
\begin{itemize}

\item
{\it Maximum Angle condition.} Let the maximum interior angle in any triangle $T\in\mathcal T$
be uniformly bounded by some positive $\alpha_0<\pi$.
\smallskip

\item
{\it Local Element Orientation condition.}
For any $z\in\mathcal N$, there is a rectangle $R_z\supset \omega_z$ such that $|R_z|\sim |\omega_z|$.
Furthermore, if $z\in\mathcal N\cap\pt\Omega$ is not a corner of $\Omega$, then $R_z$
has a side parallel to the segment ${\mathcal S}_z\cap\pt\Omega$.
\smallskip

\item
Also, let the number of triangles containing any node be uniformly bounded.

\end{itemize}
Note that the above 
conditions are automatically satisfied by
shape-regular triangulations.

Additionally,  we restrict our analysis to
the following two
node types
defined
using a fixed 
small constant $c_0$ (to distinguish between anisotropic and isotropic elements),
with the notation $a \ll b$ for $a<c_0b$.%
\smallskip

(1) {\it Anisotropic Nodes}, the set of which is denoted by ${\mathcal N}_{\rm ani}$, are such that
\beq\label{ani_node}
h_z\ll H_z,
\qquad\quad
h_T\sim h_z\;\;\mbox{and}\;\;H_T\sim H_z \quad\;\forall \ T\subset\omega_z.
\eeq
Note that the above implies that ${\mathcal S}_z$ contains at most two edges of length${}\lesssim h_z$ (see also Fig.\,\ref{fig_local}, left).%
\smallskip

(2) {\it Regular Nodes}, the set of which is denoted by ${\mathcal N}_{\rm reg}$, are those surrounded by
shape-regular mesh elements.
\smallskip

Note that most of our analysis applies to more general node types that were considered in \cite{Kopt15,Kopt_17_NM}; see Remarks~\ref{rem_gen_nodes1} and \ref{rem_gen_nodes2} for details.

\section{Sharper scaled trace theorem for anisotropic elements}\label{sec_trace}

Our task is to get an improved bound for the jump residual terms (see $I$ in \eqref{error} below).
The key to this will be to employ the following sharper version of the scaled trace theorem for anisotropic elements, which is the main result of this section.


\begin{lemma}\label{lem_scaled_new}
For any node $z\in\mathcal{N}={\mathcal N}_{\rm ani}\cup{\mathcal N}_{\rm reg}$, 
any function $v\in W^{1}_1(\omega_z)$, and any edge $S\subset\gamma_z$, one has
\begin{align}
\label{scaled_new_a}
\| v \phi_z\|_{1\,;S}&\lesssim \|\nabla v\|_{1\,;\omega_z}+\|v\|_{1\,;\omega_z}
\left\{\begin{array}{cl}
H_z^{-1}&\mbox{if~}S\subset\mathring\gamma_z\,,\\[2pt]
h_z^{-1}&\mbox{if~}S\subset\gamma_z\backslash\mathring\gamma_z\,,
\end{array}\right.\\[0.2cm]
\label{scaled_new_b}
|S|^{-1}\| v \phi_z\|^2_{1\,;S}&\lesssim
\|v\|_{2\,;\omega_z}\|\nabla v\|_{2\,;\omega_z}+\|v\|^2_{2\,;\omega_z}
\left\{\begin{array}{cl}
H_z^{-1}&\mbox{if~}S\subset\mathring\gamma_z\,,\\[2pt]
h_z^{-1}&\mbox{if~}S\subset\gamma_z\backslash\mathring\gamma_z\,,
\end{array}\right.
\end{align}
where $\phi_z$ is the hat basis function associated with $z$.
\end{lemma}

\begin{remark}\label{rem_sec_trace}
Similar versions of the scaled trace theorem for anisotropic elements were obtained in \cite[Lemma~3.1]{Kopt15} and \cite[\S3]{Kopt_17_NM}.
Lemma~\ref{lem_scaled_new} is an improvement in the sense that in the case of long edges
(i.e. $S\subset\gamma_z\backslash\mathring\gamma_z$), the weights at $\|\nabla v\|_{p\,;\omega_z}$ 
are sharper.
To be more precise, the version of \eqref{scaled_new_a} in \cite[Lemma~3.1]{Kopt15} has the weight $H_z/h_z\gg 1$ at $\|\nabla v\|_{1\,;\omega_z}$,
while the version of \eqref{scaled_new_b} given by \cite[Corollory~3.2]{Kopt_17_NM} also involves the weight $H_z/h_z\gg 1$ at $\|\nabla v\|_{2\,;\omega_z}$.
Importantly,  for the shaper bounds of Lemma~\ref{lem_scaled_new} to hold true, one needs to estimate  $\| v \phi_z\|_{p\,;S}$ rather than $\| v \|_{p\,;S}$ bounded in \cite{Kopt15,Kopt_17_NM}.
Note that this improvement is crucial for getting an improved weight in the jump residual part of our estimator.
\end{remark}

\begin{remark}\label{rem_gen_nodes1}
An inspection of the proof of Lemma~\ref{lem_scaled_new} shows that this lemma remains valid for the more general node types
introduced in \cite[\S2]{Kopt_17_NM}.
\end{remark}

To prove Lemma~\ref{lem_scaled_new}, we shall employ the following auxiliary result.
\begin{lemma}\label{lem_trace_phi}
For any sufficiently smooth function $v\ge 0$ on a triangle $T$ with vertices $z$, $z'$ and $z''$ and their respective opposite edges
$S$, $S'$ and $S''$, one has
\begin{subequations}
\begin{align}\label{aux_a}
\sin\angle(S',S'')\,\| v \phi_z\|_{1\,;S'}&\lesssim  \|\nabla v\|_{1\,;T} +|S''|^{-1}\|v\|_{1\,;T}  \,,\\[0.2cm]
\label{aux_b}
|S'|^{-1}\| v \phi_z\|_{1\,;S'}&\lesssim |S''|^{-1}\| v \phi_z\|_{1\,;S''}+ |S||T|^{-1}\|\nabla v\|_{1\,;T}\,.
\end{align}
\end{subequations}
\end{lemma}

\begin{figure}[b!]
~\hfill
\begin{tikzpicture}[scale=0.24]
\draw[fill] (4,2.4) circle [radius=0.4] node[above left]{$z$};
\path[draw]  (3.8,0)
--(15.8,3.5)--(4,2.4)--cycle;
%
\path[draw] 
(4,2.4)--(15.8,3.5);
\node[above] at (10,2.7) {$S$};
\end{tikzpicture}
\hfill
\begin{tikzpicture}[scale=0.24]
\draw[fill] (4,2.4) circle [radius=0.4] node[above left]{$z$};
\path[draw]  (3.8,0)--(16,0)--
(15.8,3.5)--(4,2.4)--cycle;
\path[draw] (4,2.4)--(16,0);
\path[draw] 
(4,2.4)--(15.8,3.5);
\node[above] at (10,2.7) {$S$};
\node[above] at (13.5,0.3) {$S''$};
\end{tikzpicture}
\hfill
\begin{tikzpicture}[scale=0.24]
\draw[fill] (4,2.4) circle [radius=0.4] node[above left]{$z$};
\path[draw]  (3.8,0)--(16,0)--(17,1.7)--(15.8,3.5)--(4,2.4)--cycle;
\path[draw] (4,2.4)--(16,0);
\path[draw] (17,1.7)--(4,2.4)--(15.8,3.5);
\node[above] at (10,2.7) {$S$};
\node[above] at (15,1.5) {$S''$};
\end{tikzpicture}
\hfill~
\caption{Illustration to the proof of \eqref{scaled_new_a} in Lemma~\ref{lem_scaled_new}:
case (i) (left); case (ii) with a single application of \eqref{aux_b} (centre);
case (ii) with a double application of \eqref{aux_b} (right).}
\label{fig_cases}
\end{figure}
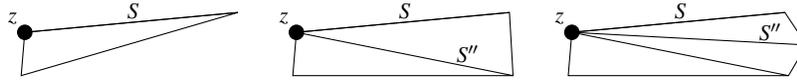

\begin{proof}
For \eqref{aux_a},
let $\bmu''$ be the unit vector along $S''$ directed from $z'$ to $z$ so that $\nabla \phi_z\cdot\bmu''=|S''|^{-1}$.
Note that $\nabla\cdot (v\phi_z \bmu'')=\nabla (v\phi_z)\cdot \bmu''$, so the divergence theorem yields
$$
\int_{\pt T}(v\phi_z\bmu'')\cdot \bnu =\int_T \nabla (v\phi_z)\cdot \bmu''=\int_T\bigl(\phi_z\nabla v\cdot \bmu'' +|S''|^{-1}v \bigr).
$$
Here, to evaluate the integral $\int_{\pt T}$, note that $\bmu''\cdot \bnu=0$ on $S''$ and $\phi_z=0$ on $S$,
while $\bmu''\cdot \bnu=\sin\angle(S',S'')$ on $S'$,
so
$\int_{\pt T}(v\phi_z\bmu'')\cdot \bnu=\sin\angle(S',S'')\int_{S'}v\phi_z $.
The desired bound \eqref{aux_a} follows.

To get \eqref{aux_b}, we modify the proof of \cite[Lemma~7.1]{Kopt_17_NM}.
Set $w=v\phi_z$ and also ${\mathcal A}_{S}\,w:=|S|^{-1}\int_S w$ for any edge $S$.
Now, with the $\zeta$-axis having the inward normal direction to $S$, and $\hbar:=2|T||S|^{-1}$, one gets
${\mathcal A}_{S'}\,w-{\mathcal A}_{S''}\,w=\hbar^{-1} \int_0^{\hbar}\! \,\bigl(w|_{S'}-w|_{S''}\bigr)\,d\zeta$.
This yields \eqref{aux_b} as $\phi_z$ does not change in the direction normal to $\zeta$.
\proofend
\end{proof}

\noindent{\it Proof of Lemma~\ref{lem_scaled_new}.}
First, note that \eqref{scaled_new_b} follows from \eqref{scaled_new_a} as
$|S|^{-1}\| v \phi_z\|^2_{1\,;S}\le \| v^2 \phi_z^2\|_{1\,;S}\le \| v^2 \phi_z\|_{1\,;S}$, while $\nabla v^2=v\nabla v$.
With regard to \eqref{scaled_new_a}, it suffices to prove it for the case $v\ge 0$, as
if $v$ changes sign on $\omega_z$, apply \eqref{scaled_new_a} with $v$ replaced by $v_\tau:=\sqrt{v^2+\tau^2}\ge 0$, where $\tau$ is a small positive constant (while $|\nabla v_\tau|\le|\nabla v|$),
 and then let $\tau\rightarrow 0^+$ so that $v_\tau\rightarrow |v|$.

Thus it remains to show \eqref{scaled_new_a} for $v\ge0$. When $S\subset\mathring{\gamma}_z$, this bound follows from a similar bound on $\| v\|_{1\,;S}$
 in \cite[Lemma~3.1]{Kopt15}. Now consider $S\subset\gamma_z\backslash\mathring{\gamma}_z$.
Then $S$ is a long edge shared by two anisotropic triangles. Consider two cases; see Fig.\,\ref{fig_cases}.
Case~(i):
If in at least one of these triangles, $T$, the angle at $z$ is $\gtrsim 1$, then an application of \eqref{aux_a} yields
$\| v \phi_z\|_{1\,;S}\lesssim \|\nabla v\|_{1\,;T}+h_z^{-1}\|v\|_{1\,;T}$, and \eqref{scaled_new_a} follows.
Case~(ii):
Otherwise, in any triangle $T$ sharing the edge $S$, the other edge $S''$ originating at $z$ is also of length $\sim H_z$, while the edge opposite to $z$ is of length $\sim h_z$.
Then
an application of \eqref{aux_b} yields
$H_z^{-1}\| v \phi_z\|_{1\,;S}\lesssim H_z^{-1}\| v \phi_z\|_{1\,;S''}+ H_z^{-1}\|\nabla v\|_{1\,;T}$
or, equivalently, $\| v \phi_z\|_{1\,;S}\lesssim \| v \phi_z\|_{1\,;S''}+ \|\nabla v\|_{1\,;\omega_z}$.
Thus, a possibly repeated application of \eqref{aux_b} reduces this case to case (i); see Fig.\,\ref{fig_cases}.
\proofend

\section{A posteriori error bounds for jump residual terms}\label{sec_error_an}
Assuming $\vvvert u_h-u \vvvert_{\eps\,;\Omega}>0$, let
\beq\label{G_def}
G:=\frac{u_h-u\ \,}{\ \vvvert u_h-u \vvvert_{\eps\,;\Omega}}
\;\;\Rightarrow \;\; \vvvert G \vvvert_{\eps\,;\Omega}=1,
\qquad\quad
g:= G-G_h,
\eeq
where $G_h\in  S_h$ is some interpolant of $G$.
Now,
a relatively standard calculation yields the following error representation \cite[\S4]{Kopt_17_NM}
\begin{align}\notag
\vvvert u_h-u \vvvert_{\eps\,;\Omega}\lesssim{}
&\!\!\sum_{z\in{\mathcal N}}\eps^2\int_{\gamma_z}(g-\bar g_z)\phi_z \llbracket\nabla u_h\rrbracket \cdot\bnu
+\!\!\sum_{z\in{\mathcal N}}\int_{\omega_z} \!\!f_h^I\,(g-\bar g_z)\phi_z+|\langle f_h-f_h^I,G\rangle|
\\
=:{}&\,I+I\!I+{\mathcal E}_{\rm quad}\,,
\label{error}
\end{align}
which holds for any $G_h\in  S_h$ and any set of real numbers $\{\bar g_z\}_{z\in\mathcal N}$ such that $\bar g_z=0$ whenever $z\in\pt\Omega$.
(To be precise, $\bar g_z$  will be specified later as a certain average of $g= G-G_h$ near $z$.)
Here $\phi_z$ denotes the standard hat basis function corresponding to $z \in \mathcal N$.

In the following proofs it will be convenient to use, with $p=1,2$, the scaled $W^{1}_p({\mathcal D})$ norm defined by
$$
\vvvertN\, v \vvvertN_{p\,;{\mathcal D}}:=\|\nabla v \|_{p\,;{\mathcal D}}+({\rm diam} {\mathcal D})^{-1}\|v \|_{p\,;{\mathcal D}}
\;\Rightarrow\;
\vvvertN\, v \vvvertN_{p\,;\omega_z}=\|\nabla v \|_{p\,;\omega_z}+H_z^{-1}\|v \|_{p\,;\omega_z}\,.
$$

\subsection{Jump residual for a partially structured anisotropic mesh}\label{sec_struct}

To illustrate our approach in a simpler setting, we first present a version of the analysis for a simpler, partially structured,
anisotropic mesh in a square domain $\Omega=(0,1)^2$.
So, throughout this section,  we make the following triangulation assumptions.%

\begin{figure}[!b]
~~~\vspace{-0.3cm}

~\hfill
\begin{tikzpicture}[scale=0.25]
\draw[ultra thick 
] (0,0) -- (21,0) ;
\path[draw,help lines]  (19,-0.2)node[below] {$\color{black}x_{i-1}$}--(19,9.5);
\path[draw,help lines]  (9.5,-0.2)node[below] {$\color{black}x_i$}--(9.5,9.5);
\path[draw,help lines]  (2,-0.2)node[below] {$\color{black}x_{i+1}$}--(2,9.5);
\path[draw]  (19,0)--(19,9.9);
\path[draw]  (9.5,0)--(9.5,9.9);
\path[draw]  (2,0)--(2,9.9);
\path[draw]  (2,1)--(9.5,1.8)--(19,0)--(9.5,0)--cycle;
\path[draw]  (9.5,1.8)--(19,1.9);
\path[draw]  (2,3)--(9.5,4.5)--(19,4)--(9.5,1.8)--cycle;
\path[draw]  (2,3)--(9.5,3)--(19,4);
\path[draw]  (2,5)--(9.5,6.4)--(19,7.5)--(9.5,4.5)--cycle;
\path[draw]  (2,6.7)--(9.5,6.4);
\path[draw]  (19,5.9)--(9.5,4.5);
\path[draw]  (2,8.1)--(9.5,8.2)--(19,7.5)--(9.5,6.4)--cycle;
%
\path[draw]  (2,9.5)--(9.5,8.2)--(19,9.4);
\end{tikzpicture}
\hfill
\begin{tikzpicture}[scale=0.21]
\path[draw,help lines] (4,-1.5)--(4,12);
\path[draw,help lines] (20,-1.5)--(20,12);
\draw[fill] (4,0) circle [radius=0.4] node[below left]{$z$};
\draw[fill] (4,5) circle [radius=0.4] node[below left]{$z'$};
\draw[fill] (20,10) circle [radius=0.4] node[above right]{$\hat z'$};
\draw[fill] (20,7) circle [radius=0.4] node[below right]{$\hat z$};
\path[draw]  (4,0)--(20,7);
\path[draw] (20,10)--(4,5);
\path[draw,<->] (21,7)--(21,10) node [midway, right] {$\gtrsim h_{\hat z}\sim h_z$};
\path[draw,<->] (3,5)--(3,10)node [midway, left] {$\lesssim h_{z'}$};
\path[draw,dashed] (3.5,7)--(22,7);
\path[draw,dotted] (2,10)--(22,10);
\path[draw,dotted] (2,5)--(4,5);

\node at (4,-2) {\color{white}.};
%
%
\end{tikzpicture}
\hfill~%
\caption{Partially structured anisotropic mesh (left);
illustration for Remark~\ref{rem_J_sturct} (right):
for any fixed edge $z\hat z$ and any edge $z'\hat z'$ intercepting the dashed horizontal line via $\hat z$,
the figure shows that $h_z\lesssim h_{z'}$, so there is a uniformly bounded number of edges of type $z'\hat z'$, so $\omega_z^*\subset \omega_z^{(J)}$ with $J\lesssim 1$.}
\label{fig_partial}
\end{figure}
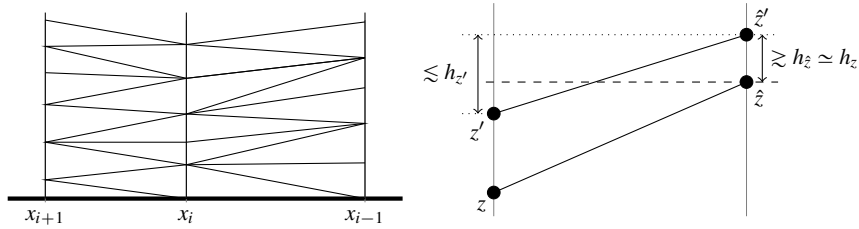

{
\begin{enumerate}
\item[A1.]
Let $\{x_i\}_{i=0}^n$ be an arbitrary mesh on the interval $(0,1)$ in the $x$ direction.
Then, let each $T\in\mathcal T$, for some $i$,\\
(i) have the shortest edge on the line $x=x_i$;\\
(ii) have a vertex on the line $x=x_{i+1}$ or $x=x_{i-1}$
(see Fig.\,\ref{fig_partial}, left).
\smallskip

\item[A2.] Let ${\mathcal N}={\mathcal N}_{\rm ani}$, i.e.
each mesh node $z$ satisfies (\ref{ani_node}).
\smallskip

\item[A3.]
{\it Quasi-non-obtuse anisotropic elements.}
Let the maximum angle in any triangle be bounded by $\frac\pi2+\alpha_1 \frac{h_T}{H_T}$
 for some positive constant $\alpha_1$.

\end{enumerate}
}

\noindent
These conditions essentially imply that all mesh elements are anisotropic and aligned in the $x$-direction.
They also imply that if $x_z=x_i$, then
\beq\label{omega_z_star}
\omega_z\subseteq \omega_z^*:=(x_{i-1},x_{i+1})\times(y_z^-,y_z^+),
\quad
y^+_z-y_z^-\sim h_z,\quad
{\rm diam}\,\omega_z^*\sim H_z\,,
\eeq
where  $(y_z^-,y_z^+)$  is the
range of $y$ within
 $\omega_z$,
while $x_{-1}:=x_0$ and $x_{n+1}:=x_n$.

\begin{remark}\label{rem_J_sturct}
The above conditions (in particular A3) imply that
there is $J\lesssim 1$ such that
$\omega_z^*\subset \omega_z^{(J)}$ for all $z\in \mathcal{N}$,
where
$\omega_z^{(0)}:=\omega_z$, and
$\omega_z^{(j+1)}$ 
denotes the patch of elements in/touching $\omega_z^{(j)}$.
This conclusion is illustrated on Fig.\,\ref{fig_partial} (right).
(Note that $J=1$ for any non-obtuse triangulation, i.e. for the case $\alpha_1=0$ in A3.)
\end{remark}

Following~\cite{Kopt15,Kopt_17_NM},
the choice of $\bar g_z$
in (\ref{error})
is related to the orientation of anisotropic elements, and is crucial in our analysis.
Let $\bar g_z=0$ for $z\in\partial \Omega$, and, otherwise, for $x_z=x_i$ with some $1\le i\le n-1$, let
\beq\label{g_bar_struct}
\int_{x_{i-1}}^{x_{i+1}}(g(x,y_z)-\bar g_z)\, \varphi_i(x)\,dx=0.
\eeq
Here we use the standard one-dimensional hat function $\varphi_i(x)$ associated with the mesh $\{x_i\}$
(i.e. it has support on $(x_{i-1},x_{i+1})$, equals $1$ at $x=x_i$, and is linear on
$(x_{i-1},x_{i})$ and $(x_{i},x_{i+1})$).

\begin{theorem_}\label{lem_J_struct}
Let
$g=G-G_h$
with $G$ from~\eqref{G_def} and any $G_h\in S_h$, while
\beq\label{Theta_A123}
\Theta:=\eps^2\|\nabla g \|^2_{2\,;\Omega}+\sum_{z\in{\mathcal N}}\bigl(1+\eps^2 H_z^{-2}\bigr) \|g \|^2_{2\,;\omega_z}\,.
\eeq
Then $\vvvert u_h-u \vvvert_{\eps\,;\Omega}\lesssim I+I\!I+{\mathcal E}_{\rm quad}$, where
the right-hand side terms are specified in \eqref{error},
and,
under conditions A1--A3,
\begin{align}\label{I_struct}
|I|&\lesssim
\Bigl\{\Theta\,\sum_{z\in{\mathcal N}}
\bigl[
\min\{|\omega_z|,\,\eps h_z \}
\,\bigl(\eps \mathring{J}_z\bigr)^2
+\min\{|\omega_z|,\,\eps H_z \}
\,\bigl(\eps J_z\bigr)^2
\bigr]
\Bigr\}^{1/2}\,,
\end{align}
 where
$\mathring{J}_z:=\|\llbracket\nabla u_h\rrbracket\|_{\infty\,;\mathring{\gamma}_z}$
and $J_z:=\|\llbracket\nabla u_h\rrbracket\|_{\infty\,;\gamma_z\backslash\mathring{\gamma}_z}\,$.
\end{theorem_}

\begin{corollary_}
Under conditions A1--A3, one has \eqref{first result} with $\lambda_z=\eps H_z\,$.
\end{corollary_}
\begin{proof}
To get the desired result, combine \eqref{I_struct} with the bound \cite[(5.8)]{Kopt_17_NM} on $I\!I$,
the straightforward bound $|{\mathcal E}_{\rm quad}|
\le \| f_h-f_h^I\|_{2\,;\Omega}\,,$
and
$\Theta\lesssim\vvvert G \vvvert_{\eps\,;\Omega}= 1$
(the latter is given by \cite[Theorem~7.4]{Kopt_17_NM} under more general conditions than A1--A3).
 \proofend
\end{proof}
\smallskip

\noindent
{\it Proof of Theorem~\ref{lem_J_struct}.}
Split $I$ of (\ref{error})
as $I=\sum_{z\in{\mathcal N}}(\mathring{I}_z+I_z)$,
where\vspace{-0.1cm}
\beq\label{I_z_new}
\mathring{I}_z:=\eps^2\int_{\mathring{\gamma}_z}(g-\bar g_z)\phi_z \llbracket\nabla u_h\rrbracket \cdot\bnu,
\qquad
I_z:=\eps^2\int_{\gamma_z\backslash\mathring{\gamma}_z}(g-\bar g_z)\phi_z \llbracket\nabla u_h\rrbracket \cdot\bnu.
\eeq

First, consider $\bar g_z$,
the definition of which \eqref{g_bar_struct} implies that
$H_z |\bar g_z|\lesssim \|g \varphi_i \|_{1\,;\bar S_z}$, where $\bar S_z$ is the segment joining the points $(x_{i-1},y_z)$ and $(x_{i+1},y_z)$, so $|\bar S_z|\sim H_z$.
Versions of \eqref{scaled_new_a} and \eqref{scaled_new_b} then respectively yield
\beq\label{two_bounds_g_z}
H_z |\bar g_z|\lesssim \|\nabla g  \|_{1\,;\omega_z^*}+h_z^{-1}\|g  \|_{1\,;\omega_z^*},\;\;
H_z |\bar g_z|^2\lesssim \| g  \|_{2\,;\omega_z^*}(\|\nabla g  \|_{2\,;\omega_z^*}+h_z^{-1}\|g  \|_{2\,;\omega_z^*}).
\eeq
These two bounds will be used when estimating  both $\mathring{I}_z$ and $I_z$.

We now proceed to estimating $\mathring{I}_z$. Note that \eqref{scaled_new_a} implies that
$\|(g-\bar g_z)\phi_z\|_{1\,;\mathring{\gamma}_z}\lesssim  \vvvertN\, g \vvvertN_{1\,;\omega_z^*}\lesssim |\omega_z|^{1/2}  \vvvertN\, g \vvvertN_{2\,;\omega_z^*} $,
where we also used $\|\bar g_z\phi_z\|_{1\,;\mathring{\gamma}_z}\sim h_z |\bar g_z|$ combined with the first bound from \eqref{two_bounds_g_z}.
Similarly, 
$\|(g-\bar g_z)\phi_z\|_{1\,;\mathring{\gamma}_z}^2\lesssim h_z\| g \|_{2\,;\omega_z^*}\vvvertN\, g \vvvertN_{2\,;\omega_z^*}$,
where we employed \eqref{scaled_new_b} and the second bound from \eqref{two_bounds_g_z}.
Now, from the definition of $\mathring{I}_z$ in \eqref{I_z_new} combined with the two bounds on $\|(g-\bar g_z)\phi_z\|_{1\,;\mathring{\gamma}_z}$, one
concludes that
$$
|\mathring{I}_z|\lesssim \mathring{\theta}_z^{1/2}\,\mathring{\lambda}_z^{1/2}\,(\eps\mathring{J}_z),
\qquad
\mathring{\theta}_z:=\frac{\eps^2\min\bigl\{|\omega_z| \vvvertN\, g \vvvertN_{2\,;\omega_z^*}^{\,2}\,,\,
h_z\| g \|_{2\,;\omega_z^*}\vvvertN\, g \vvvertN_{2\,;\omega_z^*} \bigr\}}
{\mathring{\lambda}_z}\,.
$$
Set $\mathring{\lambda}_z:=\min\{|\omega_z|,\eps h_z\}$. Then,
to get the bound of type \eqref{I_struct} for $\sum_{z\in\mathcal N} \mathring{I}_z$, it remains to show that $\sum_{z\in\mathcal N}\mathring{\theta}_z\lesssim \Theta$.
For the latter, in view of
\beq\label{ab_any}
\min\{aa',bb'\}/\min\{a',b'\}\le a+b\qquad \forall\ a,a',b,b'>0,
\eeq
one gets
$\mathring{\theta}_z\lesssim \eps^2 \vvvertN\, g \vvvertN_{2\,;\omega_z^*}^{\,2}+\eps \| g \|_{2\,;\omega_z^*}\vvvertN\, g \vvvertN_{2\,;\omega_z^*}$,
which leads to  $\sum_{z\in\mathcal N}\mathring{\theta}_z\lesssim \Theta$, also using Remark~\ref{rem_J_sturct}.

For $I_z$, first, recall the bound
$|I_z|\lesssim \eps  \vvvertN\, g \vvvertN_{1\,;\omega_z^*}(\eps J_z)$
from \cite[(5.12)]{Kopt_17_NM}, which implies  $|I_z|\lesssim \eps  |\omega_z|^{1/2}\vvvertN\, g \vvvertN_{2\,;\omega_z^*}(\eps J_z)$.
An alternative bound on $I_z$ follows from
 $\|(g-\bar g_z)\phi_z\|^2_{1\,;\gamma_z\backslash\mathring{\gamma}_z}\lesssim H_z\| g \|_{2\,;\omega_z^*}(\| \nabla g \|_{2\,;\omega_z^*}+h_z^{-1}\| g \|_{2\,;\omega_z^*})$,
where the latter is obtained by an application of \eqref{scaled_new_b} for $g$, while the second bound from \eqref{two_bounds_g_z} is employed for $\bar g_z$.
Combining the two bounds on $I_z$, we arrive at
\begin{align}
|I_z|&{}\lesssim \theta_z^{1/2}\,\lambda_z^{1/2}\,(\eps J_z),
\notag\\\label{theta_z1}
\theta_z&{}:=\frac{\eps^2\min\bigl\{|\omega_z| \vvvertN\, g \vvvertN_{2\,;\omega_z^*}^{\,2}\,,\,
H_z\| g \|_{2\,;\omega_z^*}(\| \nabla g \|_{2\,;\omega_z^*}+h_z^{-1}\| g \|_{2\,;\omega_z^*}) \bigr\}}
{\lambda_z}\,.
\end{align}
Here set $\lambda_z:=\min\{|\omega_z|,\,\eps H_z(1+\eps h_z^{-1})\}$.
Now, again using \eqref{ab_any}, one gets
\beq\label{theta_z2}
\theta_z\lesssim \eps^2 \vvvertN\, g \vvvertN_{2\,;\omega_z^*}^{\,2}
+
\eps \| g \|_{2\,;\omega_z^*}(\| \nabla g \|_{2\,;\omega_z^*}+h_z^{-1}\| g \|_{2\,;\omega_z^*})/(1+\eps h_z^{-1}),
\eeq
and hence $\sum_{z\in\mathcal N}\theta_z\lesssim \Theta$.
Finally, to get the bound of type \eqref{I_struct} for $\sum_{z\in\mathcal N} {I}_z$, it remains to note that
 $\lambda_z=\min\{|\omega_z|,\,\eps H_z[1+\eps h_z^{-1}]\}\sim \min\{|\omega_z|,\,\eps H_z\}$.
  \proofend

\begin{remark}
While the definition \eqref{g_bar_struct} for $\bar g_z$ is quite different from a
standard choice
(see, e.g., \cite[Lecture~5]{[Nochetto_lecture_notes]}),
its role
may not be immediately obvious
in the proof of Theorem~\ref{lem_J_struct}.
To clarify this, note that it is crucial for the bound $|I_z|\lesssim \eps  \vvvertN\, g \vvvertN_{1\,;\omega_z^*}(\eps J_z)$ quoted from \cite[(5.12)]{Kopt_17_NM}.
To be more precise, the latter bound is obtained in \cite{Kopt_17_NM} using the representation
\begin{align}\notag
I_z=I_z'+I_z''+I_z'''
&:=
\eps^2\int_{\gamma_z\backslash\mathring{\gamma}_z}(g-\bar g_z)\phi_z \llbracket\partial_x u_h\rrbracket\, \bnu_x
\\\notag
&{}+
\eps^2\int_{\gamma_z\backslash\mathring{\gamma}_z}[g- g(x,y_z)]\,\phi_z \llbracket\partial_y u_h\rrbracket\, \bnu_y
\\\notag
&{}+
\eps^2\int_{\gamma_z\backslash\mathring{\gamma}_z}[g(x,y_z)- \bar g_z]\,\phi_z \llbracket\partial_y u_h\rrbracket\, \bnu_y\,,
\end{align}
where
$\llbracket w\rrbracket$, for any $w$, 
is understood as the jump in $w$ across any edge in $\gamma_z$ evaluated in the {anticlockwise} direction about~$z$.
Importantly, here $I_z'''=0$ due to our choice of $\bar g_z$ (as well as due to the partial structure of our mesh; in a more general case, the estimation of $I_z'''$
is more intricate).
\end{remark}

\subsection{Jump residual for for general anisotropic meshes}\label{ssec_gen_anal}

\begin{theorem_}\label{lem_J_gen}
Suppose that $\mathcal{N}={\mathcal N}_{\rm ani}\cup{\mathcal N}_{\rm reg}$ and all corners of $\Omega$ are in ${\mathcal N}_{\rm reg}$.
Let
$g=G-G_h$
with $G$ from~\eqref{G_def} and any $G_h\in S_h$, while
$\Theta$ is defined by \eqref{Theta_A123}. 
Then $\vvvert u_h-u \vvvert_{\eps\,;\Omega}\lesssim I+I\!I+{\mathcal E}_{\rm quad}$, where
the right-hand side terms are specified in \eqref{error},
and
\begin{align}\label{I_gen_gen}
|I|&\lesssim
\Bigl\{\Theta\,\sum_{z\in{\mathcal N}}\min\{|\omega_z| ,\,\eps H_z\}
\,\bigl\|\eps \llbracket\nabla u_h\rrbracket \bigr\|^2_{\infty\,;\gamma_z}
\Bigr\}^{1/2}\,.
\quad
\end{align}
\end{theorem_}

\begin{corollary_}
Under the conditions of Theorem~\ref{lem_J_gen}, one has \eqref{first result} with $\lambda_z=\eps H_z\,$.
\end{corollary_}
\begin{proof}
To get the desired result, combine \eqref{I_gen_gen} with the bound \cite[(6.2)]{Kopt_17_NM} on $I\!I$,
the straightforward bound $|{\mathcal E}_{\rm quad}|
\le \| f_h-f_h^I\|_{2\,;\Omega}\,,$
and
$\Theta\lesssim\vvvert G \vvvert_{\eps\,;\Omega}= 1$
(the latter follows from \cite[Theorem~7.4]{Kopt_17_NM} as $\mathcal{N}={\mathcal N}_{\rm ani}\cup{\mathcal N}_{\rm reg}$).
 \proofend
\end{proof}

\begin{remark}\label{rem_gen_nodes2}
In view of Remark~\ref{rem_gen_nodes1},
an inspection of the proof of Theorem~\ref{lem_J_gen} shows that this theorem remains valid for the more general node types
introduced in \cite[\S2]{Kopt_17_NM}, and furthermore, can be combined with a shaper bound for the interior residual terms given by \cite[Theorem~6.2]{Kopt_17_NM}.
\end{remark}

\noindent{\it Proof of Theorem~\ref{lem_J_gen}.}
Split $I$ of (\ref{error})  as $I=\sum_{z\in\mathcal N}I_z$,
where $I_z$ is defined as in \eqref{I_z_new}, only with $\gamma_z\backslash\mathring{\gamma}_z$ replaced by $\gamma_z$.
It suffices to show that for some edge subset ${\mathcal S}^*\subset{\mathcal S}$ with some
quantities ${\mathcal I}_{S;z}$ associated with any $S\in {\mathcal S}_z\cap{\mathcal S}^*$
(to be specified below), one has
\begin{subequations}\label{I_zS}
\begin{align}\label{I_zS_sum}
\sum_{z\in\mathcal N}\,\sum_{S\in {\mathcal S}_z\cap{\mathcal S}^*}{\mathcal I}_{S;z}&=0,
\\[0.2cm]\label{I_zS_main}
|I_z+\sum_{S\in {\mathcal S}_z\cap{\mathcal S}^*}{\mathcal I}_{S;z}|&\lesssim \eps \vvvertN\, g \vvvertN_{1\,;\omega_z}\, \bigl\|\eps \llbracket\nabla u_h\rrbracket \bigr\|_{\infty\,;\gamma_z}
\notag\\[-0.2cm]
&\lesssim \eps |\omega_z|^{1/2}\vvvertN\, g \vvvertN_{2\,;\omega_z}\, \bigl\|\eps \llbracket\nabla u_h\rrbracket \bigr\|_{\infty\,;\gamma_z}\,,
\\[0.1cm]\label{I_zS_rough}
|I_z|+\sum_{S\in {\mathcal S}_z\cap{\mathcal S}^*}
|{\mathcal I}_{S;z}|&\lesssim  \eps
\Bigl\{ H_z\| g \|_{2\,;\omega_z}(\| \nabla g \|_{2\,;\omega_z}+h_z^{-1}\| g \|_{2\,;\omega_z})\Bigr\}^{1/2}
\,
\bigl\|\eps \llbracket\nabla u_h\rrbracket \bigr\|_{\infty\,;\gamma_z}\,.
\end{align}
\end{subequations}
Indeed, \eqref{I_zS_sum} implies that $I=\sum_{z\in\mathcal N} {I}_z=\sum_{z\in\mathcal N} (I_z+\sum_{S\in {\mathcal S}_z\cap{\mathcal S}^*}{\mathcal I}_{S;z})$,
while \eqref{I_zS_main}, \eqref{I_zS_rough} yield
\begin{align}
&|I_z+\sum_{S\in {\mathcal S}_z\cap{\mathcal S}^*}{\mathcal I}_{S;z}|\lesssim \theta_z^{1/2}\,\lambda_z^{1/2}\,\bigl\|\eps \llbracket\nabla u_h\rrbracket \bigr\|_{\infty\,;\gamma_z}\,,
\notag\\
\label{theta_z_gen}
&\theta_z{}:=\frac{\eps^2\min\bigl\{|\omega_z| \vvvertN\, g \vvvertN_{2\,;\omega_z}^{\,2}\,,\,
H_z\| g \|_{2\,;\omega_z}(\| \nabla g \|_{2\,;\omega_z}+h_z^{-1}\| g \|_{2\,;\omega_z}) \bigr\}}
{\lambda_z}\,.
\end{align}
Here set $\lambda_z:=\min\{|\omega_z|,\eps H_z(1+\eps h_z^{-1})\}$.
Then \eqref{theta_z_gen} becomes a version of \eqref{theta_z1} with $\omega_z^*$ replaced by $\omega_z$, so proceeding as in the proof of Theorem~\ref{lem_J_struct} (i.e. again employing \eqref{ab_any}),
one gets a version \eqref{theta_z2}  with $\omega_z^*$ replaced by $\omega_z$, which leads to $\sum_{z\in\mathcal N}\theta_z\lesssim \Theta$.
Now, to get the desired bound  \eqref{I_gen_gen}, it remains to note that
 $\lambda_z=\min\{|\omega_z|,\,\eps H_z(1+\eps h_z^{-1})\}\sim \min\{|\omega_z|,\,\eps H_z\}$.

So, to complete the proof, we need to establish \eqref{I_zS}.
 Relations \eqref{I_zS_sum} and \eqref{I_zS_main} immediately follow from \cite[(6.10),\,(6.11a),\,(6.11b)]{Kopt17}
 for a certain choice of $\{\bar g_z\}_{z\in\mathcal N}$,
 the edge subset ${\mathcal S}^*\subset{\mathcal S}$ and
the quantities ${\mathcal I}_{S;z}$ associated with any $S\in {\mathcal S}_z\cap{\mathcal S}^*$.
We need to recall their definitions
to prove the remaining required bound \eqref{I_zS_rough}
(which is a sharper version of \cite[((6.11c)]{Kopt17}).

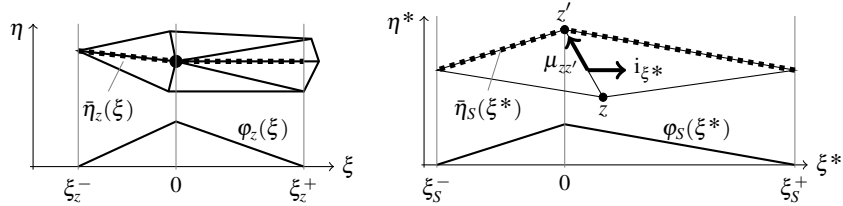
\begin{figure}[!b]
~\hfill\begin{tikzpicture}[scale=0.20]
  \draw[->] (-0.2,0) -- (20,0) node[right] {$\xi$};
  \draw[->] (0,-0.2) -- (0,9.5) ;
  \node[left] at (0,9) {$\eta$};
\draw[fill] (9.5,7) circle [radius=0.4];
\path[draw, thick]  (9,5)--(18,5)--(19,7)--(18.5,8.5)--(9.1,9)--(3,7.7)--cycle;
\path[draw, thick] (9,5)--(9.5,7)--(18,5);
\path[draw, thick] (19,7)--(9.5,7)--(18.5,8.5);
\path[draw, thick] (9.1,9)--(9.5,7)--(3,7.7);
\path[draw,dotted,
line width=2pt
] (3,7.7)--(9.5,7)--(18,7);
\path[draw,help lines]  (5,5)node[below] {$\color{black}\bar\eta_z(\xi)$}--(6,7);
\path[draw,thick] (3,0) -- (9.5,3)--(18,0);
\node[above right] at (13,1.2) {$\varphi_z(\xi)$};
\path[draw,help lines]  (18,-0.2)node[below] {$\color{black}\xi_z^+$}--(18,9.5);
\path[draw,help lines]  (9.5,-0.2)node[below] {$\color{black}0$}--(9.5,9.5);
\path[draw,help lines]  (3,-0.2)node[below] {$\color{black}\xi_z^-$}--(3,9.5);
\end{tikzpicture}
\hfill
\begin{tikzpicture}[xscale=0.17, yscale=0.18]
  \draw[->] (-14.2,-7) -- (16,-7) node[right] {$\xi^*$};
  \draw[->] (-14,-7.2) -- (-14,4) ;
  \node[left] at (-14,3.5) {$\eta^*$};
  \path[draw,thick] (-13,-7) -- (-3,-4)--(15,-7);
\node[above right] at (4,-5.8) {$\varphi_S(\xi^*)$};
\path[draw,help lines]  (15,-7.2)node[below] {$\color{black}\xi_S^+$}--(15,4);
\path[draw,help lines]  (-3,-7.2)node[below] {$\color{black}0$}--(-3,3);
\path[draw,help lines]  (-13,-7.2)node[below] {$\color{black}\xi_S^-$}--(-13,4);
\draw[fill] (0,-2) circle [radius=0.3];
\draw[fill] (-3,3) circle [radius=0.3];
%
\path[draw]  (0,-2)node[below]{$z$}--(15,0)
--(-3,3)node[above]{$z'$}--(-13,0)
--cycle;
\path[draw]  (0,-2)--(-3,3);
\path[draw,dotted,
line width=2pt
] (15,0)--(-3,3)--(-13,0);
\path[draw,help lines]  (-9,-1.5)node[below] {$\color{black}\bar\eta_S(\xi^*)$}--(-8,1.5);
\draw[line width=1.5pt, ->] (-1.2,0) -- (-2.7,2.5); \node[left] at (-1.2,0.4) {${\bmu}_{zz'}$};
\draw[line width=1.5pt, ->] (-1.2,0) -- (1.8,0) ; \node[right] at (1.8,0) {${\rm i}_{\xi^*}$};
%
%
\path[draw,white] (0,-8)--(0,-8);
\end{tikzpicture}\hfill~
\caption{Local notation associated with a node $z\in\mathcal N_{\rm ani}$ (left), and an edge $S\in{\mathcal S}^*$  with endpoints $z$ and $z'$ (right).}
\label{fig_local}
\end{figure}

First, we recall the definition of $\{\bar g_z\}_{z\in\mathcal N}$.
In view of the Local Element Orientation condition (see \S\ref{sec_nodes}),
for each fixed $z\in\mathcal N$, introduce the following {local notation}.
Let the local cartesian coordinates $(\xi,\eta)$ be such that $z=(0,0)$, and the unit vector ${\rm i}_\xi$
in the $\xi$ direction lies along the longest edge $\hat S_z\in{\mathcal S}_z$
(see Fig.\,\ref{fig_local} (left)).
For $z\in{\mathcal N}_{\rm ani}\cap\pt\Omega$ (hence $z$ is not a corner of $\Omega$),
let ${\rm i}_{\xi}$ be
either parallel
or orthogonal to $\pt\Omega$ at $z$ (depending on whether $\omega_z$ is, roughly speaking, parallel or orthogonal to $\pt\Omega$).

Next, split ${\mathcal S}_z=\mathring{\mathcal S}_z\cup{\mathcal S}_z^+\cup{\mathcal S}_z^-$,
where 
$\mathring{\mathcal S}_z=\{S\subset{\mathcal S}_z: |S|\lesssim h_z\}$
(so $\mathring{\gamma}_z=\mathring{\mathcal S}_z\setminus \pt\Omega$).
Here we also use ${\mathcal S}_z^\pm:=\{S\subset{\mathcal S}_z\setminus \mathring{\mathcal S}_z: S_\xi\subset \mathbb R_\pm\}$,
where $S_\xi={\rm proj}_\xi(S)$ denotes the projection of $S$ onto the $\xi$-axis.
Now, let
$(\xi_z^-, \xi_z^+)$ be the maximal interval such that
$(\xi_z^-, 0)\subset S_\xi$ for all $S\in{\mathcal S}_z^-$
and
$( 0,\xi_z^+)\subset S_\xi$ for all $S\in{\mathcal S}_z^+$.
Also, let
$\varphi_z (\xi)$ be the standard piecewise-linear hat-function with support on $(\xi_z^-, \xi_z^+)$
and equal to~$1$ at $\xi=0$.
Note that if ${\mathcal S}_z^-=\emptyset$ (and ${\mathcal S}_z^+=\emptyset$), then we set $\xi_z^-=0$ (and $\xi_z^+=0$) and do not use $\varphi_z$ for $\xi<0$ (and $\xi>0$).

Next, for $\xi\in[\xi_z^-,\xi_z^+]$ define a continuous function $\bar\eta_z(\xi)$ as follows:
(i)~$\bar\eta_z(\xi)$ is linear on $[\xi_z^-,0]$ and $[0,\xi_z^+]$;
(ii)~$\bar\eta_z(0)=0$;
(iii)~$(\xi,\bar\eta_z(\xi))\in\omega_z$ for all $\xi\in(\xi_z^-,\xi_z^+)$.
(For example, 
one may choose $\bar\eta_z(\xi)$ so that
$\{(\xi,\bar\eta_z(\xi)) : \xi\in(\xi_z^-,0)\}$ lies on any edge in ${\mathcal S}_z^-$,
while $\{(\xi,\bar\eta_z(\xi)) : \xi\in(0,\xi_z^+)\}$ lies on any edge in ${\mathcal S}_z^+$;
see Fig.\,\ref{fig_local} (left).)

We are now prepared to {specify $\bar g_z$}.
Let $\bar g_z:=0$ if $z\in\partial \Omega$ or $z\in{\mathcal N}_{\rm reg}$ (as for the latter, $\xi_z^-=\xi_z^+=0$), and, otherwise, let
\beq\label{g_bar_gen}
\int_{\xi_z^-}^{\xi_z^+}\bigl[g(\xi,\bar\eta_z(\xi))-\bar g_z\bigr]\, \varphi_z(\xi)\,d\xi=0.
\vspace{-0.1cm}
\eeq
%
Also, let $\bar S_z^-:=\{(\xi,\bar\eta_z(\xi)) : \xi\in(\xi_z^-,0)\}$ and
$\bar S_z^+:=\{(\xi,\bar\eta_z(\xi)) : \xi\in(0,\xi_z^+)\}$,
i.e.
$\bar S_z^\pm$ is the segment joining $(0,0)$ and $(\xi_z^\pm,\bar\eta_z(\xi_z^\pm))$.

We can now proceed to getting a bound of type \eqref{I_zS_rough} for $|I_z|$.
First, consider $\bar g_z$,
the definition of which \eqref{g_bar_gen} implies that
$H_z |\bar g_z|\lesssim \|g \varphi_z \|_{1\,;\bar S^-_z\cup \bar S^-_z}$, where $|S^-_z\cup \bar S^-_z|\sim H_z$.
Using \eqref{scaled_new_a} and \eqref{scaled_new_b} then  yields a version of
\eqref{two_bounds_g_z}, only with $\omega_z^*$ replaced by $\omega_z$
(as now $\bar S_z^-\cup\bar S_z^+\subset\omega_z$).
Next, we get
$\|(g-\bar g_z)\phi_z\|^2_{1\,;\gamma_z}\lesssim H_z\| g \|_{2\,;\omega_z}(\| \nabla g \|_{2\,;\omega_z}+h_z^{-1}\| g \|_{2\,;\omega_z})$,
which is obtained by an application of \eqref{scaled_new_b} for $g$, while the second bound from \eqref{two_bounds_g_z} is employed for $\bar g_z$.
Combining this with the definition of $I_z$ immediately yields a bound of type \eqref{I_zS_rough} for $|I_z|$.

To establish a bound of type \eqref{I_zS_rough} for $|{\mathcal I}_{S;z}|$,
we now recall the definitions of the edge subset ${\mathcal S}^*\subset{\mathcal S}$ and
the quantities ${\mathcal I}_{S;z}$ associated with any $S\in {\mathcal S}_z\cap{\mathcal S}^*$ from \cite{Kopt_17_NM}.
Let
${\mathcal S}^*:=\cup_{z\in {\mathcal N}_{\rm ani}\setminus \pt\Omega}\,\mathring{\mathcal S}_z$, and for any $S\in {\mathcal S}^*$ with endpoints $z$ and $z'$,
define
\beq\label{I_Sz}
{\mathcal I}_{S;z}:=\eps^2\alpha_S\, \bmu_{zz'}\cdot {\rm i}_{\xi^*}\, J_S\,,\quad\;
\alpha_S:=\int^{\xi_{S}^+}_{0}[g(\xi^*,\bar\eta_S(\xi^*))- \bar g_S]\,\varphi_S(\xi^*)\,d\xi^*.
\eeq
Here 
$J_S$ is the standard signed version of $|\llbracket\nabla u_h\rrbracket|$ on $S$,
$\bmu_{zz'}$ is the unit vector directed from $z$ to $z'$, and ${\rm i}_{\xi^*}$ is the unit vector along the $\xi^*$-axis. The local cartesian coordinates
$(\xi^*,\eta^*)$ are associated with $S$ and coincide with the local coordinates $(\xi,\eta)$ associated with either $z\in{\mathcal N}_{\rm ani}\setminus \pt\Omega$ or $z'\in{\mathcal N}_{\rm ani}\setminus \pt\Omega$
(at least one of them is always  in ${\mathcal N}_{\rm ani}\setminus \pt\Omega$).
The above $\alpha_S$ is defined by a version
$\int^{\xi_{S}^+}_{\xi_{S}^-}[g(\xi^*,\bar\eta_S(\xi^*))- \bar g_S]\,\varphi_S(\xi^*)\,d\xi^*=0$
of \eqref{g_bar_gen}.
The one-dimensional hat function $\varphi_S(\xi^*)$ is associated with the interval $(\xi_{S}^-,\xi_{S}^+)$; the latter is
the projection of $\omega_z\cap\omega_{z'}$ (which includes at most two triangles)
onto the $\xi^*$-axis.
The piecewise-linear function $\bar\eta_S(\xi^*)$ is defined similarly to $\bar\eta_z(\xi)$ under the restriction that
any point $(\xi^*,\bar\eta(\xi^*))\in\omega_z\cap\omega_{z'}$
(see Fig.\,\ref{fig_local}(right)).
%

Under this definition,
a bound of type \eqref{I_zS_rough} for $|{\mathcal I}_{S;z}|$ is established similarly to a similar bound for $|I_z|$.
(Note also that  $\bmu_{zz'}+\bmu_{z'z}=0$ 
in \eqref{I_Sz}, so ${\mathcal I}_{S;z}+{\mathcal I}_{S;z'}=0$, which implies
\eqref{I_zS_sum}.)
This completes the proof of \eqref{I_zS}, and hence of \eqref{I_gen_gen}.
 \proofend



\begin{thebibliography}{20}
\bibitem{AinsOd_2000}
Ainsworth,~M., Oden,~J.~T.:
A posteriori error estimation in finite element analysis.  Wiley-Interscience, New York (2000)



\bibitem{ChKopt}
Chadha,~N.~M., Kopteva,~N.:
Maximum norm a posteriori error estimate for a 3d singularly perturbed semilinear reaction-diffusion problem.
Adv. Comput. Math. {\bf35}, 33--55 (2011)

\bibitem{DK14}
Demlow,~A., Kopteva,~N.:  Maximum-norm a posteriori error estimates for singularly perturbed elliptic reaction-diffusion problems.
%
Numer. Math. {\bf133}, 707--742 (2016)



  \bibitem{Kopt_MC_07}
Kopteva,~N.: Maximum norm error analysis of a 2d singularly perturbed semilinear reaction-diffusion problem. Math. Comp. {\bf76}, 631--646 (2007)

\bibitem{Kopt08}
Kopteva,~N.:  Maximum norm a posteriori error estimate for a 2d singularly
  perturbed reaction-diffusion problem.
  SIAM J. Numer. Anal. {\bf46}, 1602--1618  (2008)



 \bibitem{Kopt15}
Kopteva,~N.: Maximum-norm a posteriori error estimates for singularly perturbed reaction-diffusion problems on anisotropic meshes.
 SIAM J. Numer. Anal. {\bf53},  2519--2544 (2015)

  \bibitem{Kopt_17_NM}
Kopteva,~N.: Energy-norm a posteriori error estimates for singularly perturbed reaction-diffusion problems on anisotropic meshes. Numer. Math., {\bf137}, 607--642 (2017)

 \bibitem{Kopt17}
Kopteva,~N.: Fully computable a posteriori error estimator using anisotropic flux equilibration on
anisotropic meshes. Submitted for publication (2017),
 arXiv:1704.04404.


\bibitem{Kopt_EOR}
Kopteva,~N., ~O'Riordan,~E.: Shishkin meshes in the numerical solution of singularly perturbed differential equations.
Int. J. Numer. Anal. Model. {\bf7}, 393--415 (2010)

\bibitem{Kunert2000}
 Kunert,~G.:
 An a posteriori residual error estimator for the finite element method on anisotropic tetrahedral meshes.
 Numer. Math.  {\bf86}, 471--490 (2000)

\bibitem{Kunert2001}
 Kunert,~G.:
 Robust a posteriori error estimation for a singularly perturbed reaction-diffusion equation on anisotropic tetrahedral meshes.
 Adv. Comput. Math.  {\bf15}, 237--259 (2001)

\bibitem{KunVer00}
 Kunert,~G., Verf{\"u}rth,~R.:
  Edge residuals dominate a posteriori error estimates
  for linear finite element methods on anisotropic triangular and tetrahedral meshes.
Numer. Math. {\bf86}, 283--303 (2000)


 \bibitem{[Nochetto_lecture_notes]}
\mbox{Nochetto,~R.H.:} Pointwise
  a posteriori error estimates for monotone semi-\-linear equations.
  Lecture Notes at 2006 CNA Summer School
Proba\-bi\-listic and Analytical Perspectives on Contemporary PDEs,
{http://www.math.cmu.edu/cna/Summer06/lecturenotes/nochetto/}

  \bibitem{RStTob}
 Roos,~H.-G., Stynes,~M., Tobiska,~T.: Robust Numerical
Methods for Singularly Perturbed Differential
          Equations. Springer,
Berlin (2008)


\bibitem{Siebert96}
Siebert,~K.~G.: An a posteriori error estimator for anisotropic refinement. Numer. Math.  {\bf73}, 373--398 (1996)


\bibitem{Ver98c}
Verf{\"u}rth,~R.:  Robust a posteriori error estimators for a
  singularly perturbed reaction-diffusion equation.
  Numer. Math. {\bf78}, 479--493
  (1998)

\end{thebibliography}
\end{document}